\documentclass[12pt]{article}
\usepackage{epsfig}
\usepackage{amsfonts}
\usepackage{amsthm}
\newtheorem{theorem}{Theorem}
\newtheorem{lemma}[theorem]{Lemma}

\newcounter{lth}
\begin{document}
\title{Cyclic colorings of plane graphs with independent faces\thanks{This research was supported by the Czech-Slovenian bilateral project MEB 090805 (on the Czech side) and BI-CZ/08-09-005 (on the Slovenian side).}}
\author{Jernej Azarija\thanks{Department of Mathematics, University of Ljubljana, Jadranska~19, 1111 Ljubljana, Slovenia. E-mails: {\tt jernej.azarija@gmail.com}, {\tt rok.erman@fmf.uni-lj.si} and {\tt matjaz.krnc@gmail.com}.}\and
        Rok Erman$^{\fnsymbol{lth}}$\and
        Daniel Kr\'al'\thanks{Institute for Theoretical Computer Science (ITI), Faculty of Mathematics and Physics, Charles University, Malostransk\'e n\'am\v est\'\i{} 25, 118 00 Prague 1, Czech Republic. E-mail: {\tt kral@kam.mff.cuni.cz}. Institute for Theoretical computer science is supported as project 1M0545 by Czech Ministry of Education.}\and
	Matja{\v z} Krnc$^{\fnsymbol{lth}}$\and
        Ladislav Stacho\thanks{Department of Mathematics, Simon Fraser University, 8888 University Dr, Burnaby, BC, V5A 1S6, Canada. E-mail: {\tt lstacho@sfu.ca}. This research was supported by NSERC grant 611368.}}
\date{}
\maketitle
\begin{abstract}
Let $G$ be a plane graph with maximum face size $\Delta^*$. If all faces
of $G$ with size four or more are vertex disjoint, then $G$ has a cyclic
coloring with $\Delta^*+1$ colors, i.e., a coloring such that all vertices
incident with the same face receive distinct colors.
\end{abstract}

\section{Introduction}
\label{sect-intro}

In 1965, Ringel~\cite{bib-ringel65} introduced the notion of
{\em $1$-planar graphs}. These are graphs that can be drawn
in the plane such that every edge is crossed by at most one
other edge. Ringel~\cite{bib-ringel65} proved that $1$-planar graphs
are $7$-colorable and conjectured that they are $6$-colorable.
Ringel's conjecture was shown to be true
by Borodin~\cite{bib-borodin84,bib-borodin95} in the 1980's.

Ringel's problem fits a framework of {\em cyclic colorings},
vertex colorings of embedded graphs such that any two vertices
incident with the same face receive distinct colors. It is easy
to see that every edge-maximal $1$-planar graph can be obtained
from a plane graph with faces of size three and four by inserting 
a pair of crossing edges into each face of size four. In the other
direction, removing pairs of crossing edges in an edge-maximal
$1$-planar graph yields a plane graph with faces of size three and
four. Hence, Borodin's result~\cite{bib-borodin84,bib-borodin95}
asserts that every plane graph with maximum face size four has
a cyclic coloring using at most six colors.

Borodin's result had been conjectured as one of the cases ($\Delta^*=4$)
in the Cyclic Coloring Conjecture of Ore and Plummer~\cite{bib-ore69+}.
The conjecture asserts that every plane graph with maximum face size $\Delta^*$
has a cyclic coloring with $\lfloor 3\Delta^*/2\rfloor$ colors.
The statement of the conjecture for $\Delta^*=3$ is equivalent
to the Four Color Theorem, proved in~\cite{bib-appel76+,bib-robertson97+}.
For $\Delta^*\ge 5$,
the best known bound of $\left\lceil 5\Delta^*/3\right\rceil$
has been obtained by Sanders and Zhao~\cite{bib-sanders01+}
improving earlier bounds of Borodin~\cite{bib-borodin92,bib-borodin99+}.
A major evidence that the conjecture is true is a recent breakthrough of
Amini, Esperet and van den Heuvel~\cite{bib-amini08+}
which extends an approach of
Havet, van den Heuvel, McDiarmid and Reed~\cite{bib-havet07+,bib-havet08+};
Amini et al.~\cite{bib-amini08+} showed that
the Cyclic Coloring Conjecture is asymptotically true, i.e.,
for every $\varepsilon>0$,
there exists $\Delta_\varepsilon$ such that every plane graph of maximum face size
$\Delta^*\ge\Delta_\varepsilon$ admits a cyclic coloring with at most
$\left(\frac{3}{2}+\varepsilon\right)\Delta^*$ colors.

The Cyclic Coloring Conjecture stipulated a lot of research,
in particular, several restrictions and generalizations of
the conjecture has been considered.
Plummer and Toft~\cite{bib-plummer87+} conjectured that the asserted
bound can be improved for $3$-connected plane graphs to $\Delta^*+2$.
The conjecture of Plummer and Toft
is known~\cite{bib-enomoto01+,bib-hornak99+,bib-hornak00+,bib-hornak+} to be true for $\Delta^*\in\{3,4\}$ and $\Delta^*\ge 18$.
In another direction, a possible generalization avoiding
the restriction of face sizes, the Facial Coloring Conjecture,
was proposed in~\cite{bib-kral05+}.
This generalization asserts that vertices of every plane
graph can be colored with at most $3\ell+1$ colors in such a way that
every two vertices joined by a facial walk of length at most $\ell$
receive distinct colors. Partial results towards proving this conjecture,
which implies the Cyclic Coloring Conjecture for odd values of $\Delta^*$,
can be found in~\cite{bib-havet+,bib-havet++,bib-kral05+,bib-kral07+}.

In this paper, we consider a different restriction of the Cyclic Coloring
Conjecture which is also motivated by colorings of graphs drawn
in the plane with restricted structure of crossings, originally
introduced by Albertson~\cite{bib-albertson08}.
Two distinct crossings are independent if the end-vertices of
every pair of crossing edges are mutually different. In particular,
if all crossings are independent, then each edge is crossed
by at most one edge, i.e., graphs with mutually independent crossings
are $1$-plane graphs.
Albertson conjectured~\cite{bib-albertson08,bib-albertson} that
every graph that can be drawn in the plane with all its crossings
independent is $5$-colorable and provided partial results towards the proof
of his conjecture (other partial results
can be found in~\cite{bib-harman,bib-wenger}).
In the cyclic coloring setting, Albertson's conjecture says that
every plane graph with faces of size three and four such
that all faces of size four are vertex-disjoint is $5$-colorable.

Albertson's conjecture has been verified by two of the authors
in~\cite{bib-kral+stacho}. A natural question is what is the least number
of colors needed if the maximum face size $\Delta^*$ is larger than four and
the faces of size four or more are still vertex disjoint.
The wheels are plane graphs of this type and thus the number of colors
needed is at least $\Delta^*+1$. We prove that this number also suffices.

\section{Overview}

Let us first introduce some additional notation. A vertex of degree $d$
is a {\em $d$-vertex} and a face incident with $k$ vertices
is a {\em $k$-face}.
The graphs we consider throughout the proof have no loops and no $2$-faces
but they can have parallel edges, in which case the degree of the vertex
is considered to be the number of edges
incident with it, not the number of its neighbors.
Two vertices are {\em cyclic neighbors} if they are incident with the same face.
The {\em cyclic degree} of a vertex $v$ is the number of distinct cyclic neighbors of $v$.

A plane graph $G$ is {\em $D$-minimal} if it has no cyclic coloring
with at most $D+1$ colors, it has maximum face size at most $D$,
all its faces of size four or more are vertex-disjoint, and $G$
has the minimal number of vertices subject to the previous constraints.
Clearly, a $D$-minimal graph is $2$-connected and
has no separating cycles of length two or
three. We will use these facts implicitly throughout the paper. 
\paragraph{}
Our goal is to show that
there is no $D$-minimal graph with $D\ge 5$ (see Theorem~\ref{thm-our}).
This will combine with the previous results to the following:

\begin{theorem}
\label{thm-main}
Every plane graph with maximum face size $\Delta^*$ whose all faces
of size four or more are vertex-disjoint has a cyclic coloring
with at most $\Delta^*+1$ colors.
\end{theorem}

The general structure of the proof is the following.
We first identify configurations that cannot appear in a $D$-minimal
graph; these configurations will be called {\em reducible configurations}.
Using the knowledge of reducible configurations, we exclude the existence
of a $D$-minimal graph by assigning each vertex and face charge in such
a way that the total amount of charge is negative. The assigned charge is then
redistributed using rules preserving its amount. The original
amount of total charge will be $-12$ and we will be able to show that
the final amount of charge of all vertices and faces is non-negative.
This will exclude the existence of a $D$-minimal graph.

\section{Reducible configurations}

In this section, we study configurations that cannot appear in a $D$-minimal graph $G$.
Let us start with a simple observation on the minimum degree of a $D$-minimal graph.

\begin{lemma}
\label{lm-mindeg}
The minimum degree of every $D$-minimal graph $G$, $D\ge 5$, is at least four.
\end{lemma}

\begin{proof}
It is straightforward to show that $G$ has no $1$-vertex. Assume that $G$
has a $d$-vertex $v$, $d\in\{2,3\}$. If $v$ is incident with $3$-faces only,
then proceed as follows: remove $v$ from $G$ and
consider a cyclic $(D+1)$-coloring of the resulting graph which exists 
by the minimality of $G$. This coloring can be extended to $v$
since the cyclic degree of $v$ is at most $3\le D$.
Hence, we assume that $v$ is incident with an $\ell$-face,
$\ell\ge 4$.

Let $w$ and $w'$ be the neighbors of $v$ incident with the $\ell$-face and
$G'$ the graph obtained from $G$ by removing $v$ and adding the edge $ww'$
if the degree of $v$ is three. Observe that the maximum face size of $G'$
does not exceed the maximum face size of $G$ and the faces of size four and
more are still vertex-disjoint.

Consider a cyclic $(D+1)$-coloring of $G'$
which exists by the minimality of $G$.
We now construct a cyclic $(D+1)$-coloring of $G$.
The vertices of $G$ distinct from $v$
preserve their colors. There are at most $D$ colors that cannot be assigned
to $v$: the colors of the $\ell-1\le D-1$ colors incident with the $\ell$-face and
the color of the third neighbor of $v$ if $v$ is a $3$-vertex. We conclude
that there is a color that can be assigned to $v$ and thus the coloring
can be completed to a cyclic $(D+1)$-coloring of $G$.
\end{proof}

In the next lemma, we look at vertices of degree four and five in $D$-minimal graphs.

\begin{lemma}
\label{lm-deg45}
Let $G$ be a $D$-minimal graph, $D \geq 5$.
Every $d$-vertex, $d\in \left\{4,5\right\}$ is incident with an $\ell$-face, $\ell \geq 4$.
\end{lemma}

\begin{proof}
Consider a $d$-vertex $v$, $d\in \left\{4,5\right\}$, contained only in
$3$-faces. Let $G'$ be the graph obtained from $G$ by removing $v$ and
triangulating the new $d$-face. By the minimality of $G$,
the graph $G'$ has a cyclic $(D+1)$-coloring.
We now extend this coloring to $G$. The vertices distinct from $v$ keep their colors.
Since the cyclic degree of $v$ in $G$ is at most $d\le D$, the coloring
can be extended to $v$ which contradicts our assumption that $G$ is $D$-minimal.
\end{proof}


Next, we show that $4$-vertices can be incident with $3$-faces and $\ell$-faces, $\ell\ge 5$, only.

\begin{lemma}
\label{lm-deg4}
Let $G$ be a $D$-minimal graph, $D \geq 5$.
No $4$-face of $G$ contains a $4$-vertex. 
\end{lemma}

\begin{figure}
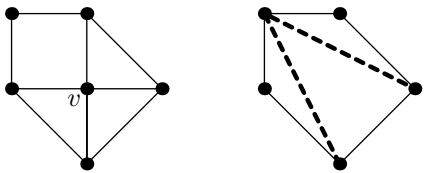

\begin{center}
\epsfbox{slika.10}
\hskip 1cm
\epsfbox{slika.11}
\end{center}
\caption{A reduction of a $4$-vertex incident with a $4$-face.}
\label{no4on4}
\end{figure}

\begin{proof}
Let $v$ be a $4$-vertex incident with a $4$-face $f$ and let $v'$ be the vertex
of the $4$-face not adjacent to $v$. By removing $v$ from $G$ and
triangulating the resulting $5$-face with edges incident with $v'$,
we obtain a graph $G'$ (see Figure \ref{no4on4}).
By the minimality of $G$, the constructed graph $G'$ has a cyclic
$(D+1)$-coloring. Since the cyclic degree of $v$ is $5$ and $D\geq 5$,
there is a color that can be assigned to $v$. This completes the coloring
to a cyclic $(D+1)$-coloring of $G$.
\end{proof}

Our next goal is to exclude the cases that a $4$-face or a $5$-face is incident
with too many $5$-vertices. This is done in the next two lemmas.

\begin{lemma}
\label{lm-4face-555}
Let $G$ be a $D$-minimal graph, $D \geq 5$.
No $4$-face of $G$ contains three $5$-vertices.
\end{lemma}

\begin{figure}
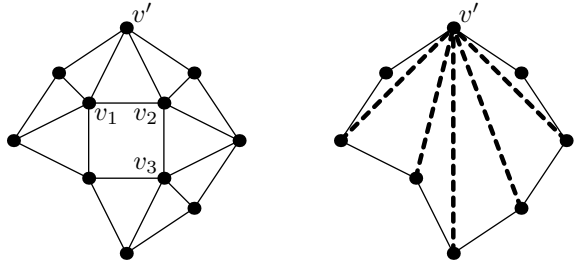

\begin{center}
\epsfbox{slika.12}
\hskip 1cm
\epsfbox{slika.13}
\end{center}
\caption{A reduction of a $4$-face incident with three $5$-vertices.}
\label{THREE5_VERTICES_1}
\end{figure}

\begin{proof}
Assume that $G$ contains a $4$-face incident with three $5$-vertices
$v_1$, $v_2$ and $v_3$ (in this order on the boundary). Let $v'$ be the
common neighbor of $v_1$ and $v_2$ (see Figure~\ref{THREE5_VERTICES_1}).
Remove the vertices $v_1$, $v_2$ and $v_3$ from $G$ and triangulate
the new $8$-face with edges originating from $v'$. Note that the obtained
graph $G'$ has no loops since $G$ has no separating triangles.
By the minimality of $G$, the graph $G'$ has a cyclic $(D+1)$-coloring.

We now extend this coloring to $G$. Since $D \geq 5$, there are
at least $6$ colors in total which can be used in the coloring.
Let $a$ be the color of $v'$. Color the vertex $v_3$
with $a$. We next color the vertices $v_1$ and $v_2$.
Each of these two vertices
has cyclic degree $6$ but two of its cyclic neighbors ($v_3$ and $v'$)
have the same color.
As $D \geq 5$, the coloring can be extended to a cyclic $(D+1)$-coloring.
\end{proof}

\begin{lemma}
\label{lm-5face}
Let $G$ be a $D$-minimal graph, $D \geq 5$.
No $5$-face $f$ of $G$ contains a vertex of degree five adjacent
to a vertex of degree four or five.
\end{lemma}

\begin{figure}
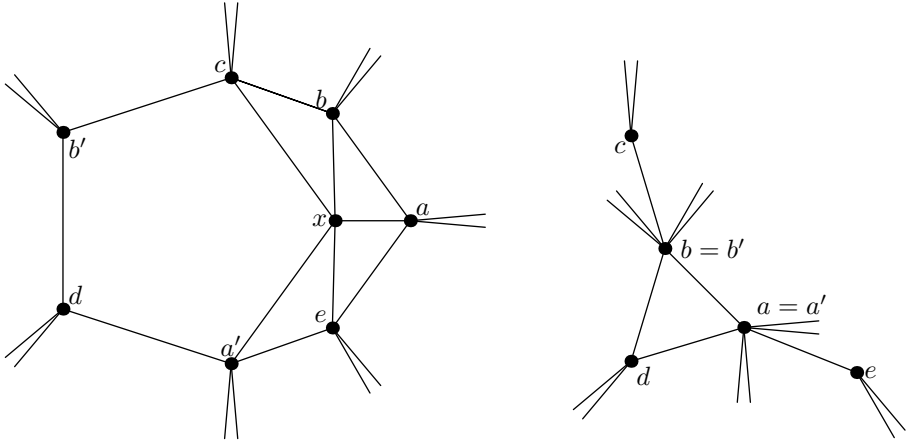

\begin{center}
\epsfbox{slika.9}
\hskip 1cm
\epsfbox{slika.16}
\end{center}
\caption{A $5$-face incident with two consecutive $5$-vertices and its reduction.}
\label{CONF_5FACE_5DEG_1}
\end{figure}

\begin{proof}
Let $x$, $c$, $b'$, $d$ and $a'$ be the vertices of $f$ (in this order) and
assume that the vertex $x$ is a $5$-vertex and $c$ is a $4$-vertex or a $5$-vertex (see Figure~\ref{CONF_5FACE_5DEG_1}).
Let $b$ be the common neighbor of $x$ and $c$, $e$ the common neighbor of $x$ and $a'$ and
$a$ the remaining neighbor of $x$. We now modify the graph $G$ to another graph $G'$.
Remove the vertex $x$, identify the vertices $a$ and $a'$, and $b$ and $b'$.
The resulting graph is $G'$ and is loopless since $G$ has no separating
triangles. By the minimality of $G$, $G'$ has a cyclic $(D+1)$-coloring.

Before extending the coloring of $G'$ to $G$, we might have to recolor
the vertex $c$ (its color can coincide with the color of the vertex $a'$ or
the color of the vertex $d$). As the cyclic degree of $c$ is at most $7\le D+2$,
it has an uncolored neighbor (the vertex $x$) and two cyclic neighbors
with the same color (the vertices $b$ and $b'$), it is possible to recolor
it.
Finally, since the cyclic degree of $x$ is $7\le D+2$ and
two pairs of its cyclic neighbors (the vertices $a$ and $a'$, and $b$ and $b'$)
have the same color, the coloring can be extended to $x$. The existence
of a cyclic $(D+1)$-coloring of $G$ contradicts the minimality of $G$.
\end{proof}

In the remaining three lemmas, we consider degrees of consecutive vertices
on an $\ell$-face, $\ell\ge 5$. We first exclude the existence of a face
with two consecutive $4$-vertices.

\begin{lemma}
\label{lm-face44}
Let $G$ be a $D$-minimal graph, $D\geq5$.
No $\ell$-face $f$, $\ell \geq 5$, of $G$ contains two consecutive $4$-vertices.  
\end{lemma}

\begin{figure}
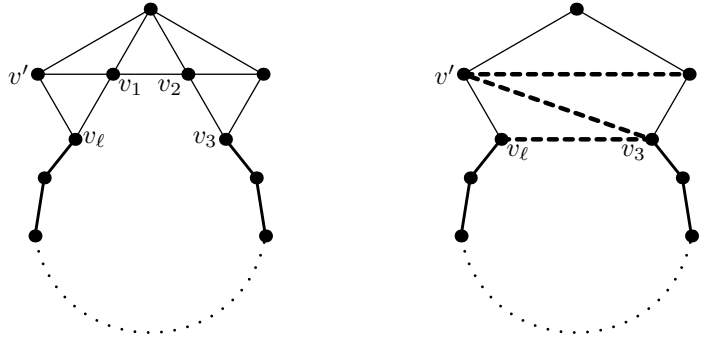

\begin{center}
\epsfbox{slika.1}
\hskip 2cm
\epsfbox{slika.2}
\end{center}
\caption{A reduction of a face with two consecutive $4$-vertices.}
\label{FIG_4VER4FACE_1}
\end{figure}

\begin{proof}
Let $v_1,\ldots,v_{\ell}$ be the vertices incident with $f$ listed in the order
on its boundary and assume that $v_1$ and $v_2$ are $4$-vertices (see Figure~\ref{FIG_4VER4FACE_1}).
Further, let $v'$ be the common neighbor of $v_1$ and $v_{\ell}$.
Form a graph $G'$ by removing $v_1$ and $v_2$, adding
the edge $v_3v_{\ell}$ and triangulating the new $5$-face by adding edges
originating from $v'$ as in Figure~\ref{FIG_4VER4FACE_1}.
Since $G$ has no separating triangles, $G'$ is loopless. Consequently,
the minimality of $G$ implies that $G'$ has a cyclic $(D+1)$-coloring.

Let $a$ be the color assigned to the vertex $v'$. If the color $a$ is assigned
to none of the vertices $v_3,\ldots,v_{\ell}$, color $v_2$ with $a$.
Otherwise color, $v_2$ with any available color (as the cyclic degree of $v_2$
is $\ell+1\le D+1$ and $v_1$ has no color, there is a color that can be used).
Observe that two cyclic neighbors of $v_1$ now have the color $a$. Since the cyclic
degree of $v_1$ is $\ell+1\le D+1$ and two of its cyclic neighbors have the same
color, the coloring can be completed to a cyclic $(D+1)$-coloring of $G$.
\end{proof}

In the final two lemmas of this section,
we exclude that one of three consecutive vertices on an $\ell$-face, $\ell\ge 5$,
would have degree four and the remaining two would have degree four or five.

\begin{lemma}
\label{lm-face45x}
Let $G$ be a $D$-minimal graph, $D\geq5$.
No $\ell$-face $f$, $\ell \geq 5$, of $G$ contains three consecutive vertices
with degrees $4,5,4$ or $4,5,5$ (in this order).
\end{lemma}

\begin{figure}
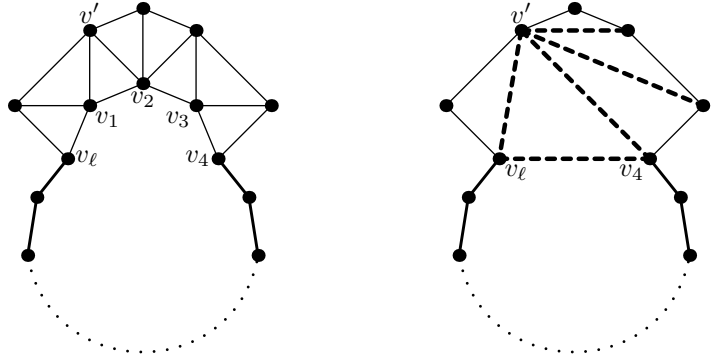

\begin{center}
\epsfbox{slika.3}
\hskip 2cm
\epsfbox{slika.4}
\end{center}
\caption{An $\ell$-face, $\ell\ge 5$, with three consecutive vertices with degrees $4,5,4$ and its reduction.}
\label{CONF_454_1}
\end{figure}

\begin{figure}
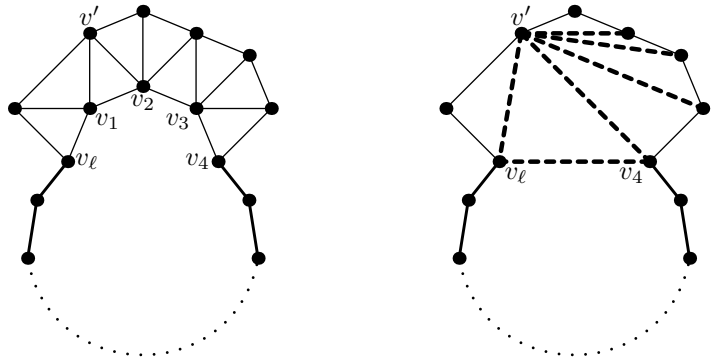

\begin{center}
\epsfbox{slika.5}
\hskip 2cm
\epsfbox{slika.6}
\end{center}
\caption{An $\ell$-face, $\ell\ge 5$, with three consecutive vertices with degrees $4,5,5$ and its reduction.}
\label{CONF_455_1}
\end{figure}

\begin{proof}
Let $v_1,\ldots,v_{\ell}$ be the vertices incident with $f$ listed in the order
on its boundary and assume that $v_1$ is a $4$-vertex, $v_2$ is a $5$-vertex and
$v_3$ is a $d$-vertex, $d\in\{4,5\}$ (see Figures~\ref{CONF_454_1} and~\ref{CONF_455_1}).
Let $G'$ be the graph obtained by removing the vertices $v_1$, $v_2$ and $v_3$,
adding the edge $v_4v_{\ell}$ and triangulating the new face by adding edges
originating from $v'$ (see the figures) where $v'$ is the common neighbor
of $v_1$ and $v_2$. Again, $G'$ has no loops as $G$ has
no separating cycles of length at most three, and
the minimality of $G$ implies that $G'$ is cyclically
$(D+1)$-colorable.

We extend a cyclic $(D+1)$-coloring of $G'$ to $G$. Let $a$ be the color assigned
to the vertex $v'$.  If the color $a$ is not assigned to any of the vertices
$v_4,\ldots,v_{\ell}$, assign $a$ to $v_3$. Otherwise, color
$v_3$ with any available color (as the cyclic degree of $v$ is at most $D+2$ and
two of its cyclic neighbors are uncolored, there is such an available color).
Color now the vertex $v_2$: the cyclic degree of $v_2$ is $D+2$
but two of its cyclic neighbors have the same color (the color $a$) and
one of its cyclic neighbors (the vertex $v_1$) is uncolored.
Finally, we color the vertex $v_1$: since its cyclic degree is $D+1$ and
two of its cyclic neighbors have the same color, there is a color that
can be assigned to $v_1$. The existence of a cyclic $(D+1)$-coloring of $G$
contradicts the minimality of $G$.
\end{proof}

\begin{lemma} 
\label{lm-face545}
Let $G$ be a $D$-minimal graph, $D\geq5$.
No $\ell$-face $f$, $\ell \geq 5$, of $G$ contains three consecutive vertices
with degrees $5,4,5$ (in this order).
\end{lemma}

\begin{figure}
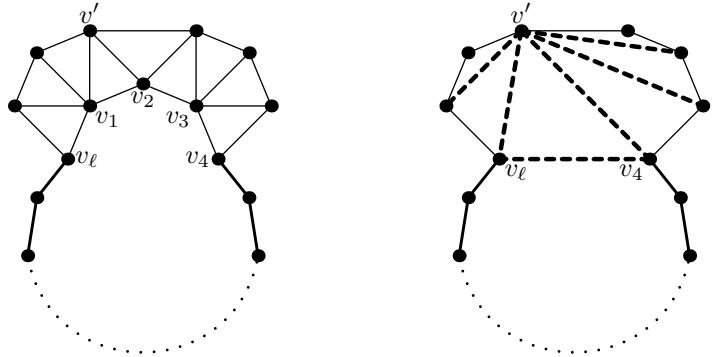

\begin{center}
\epsfbox{slika.7}
\hskip 2cm
\epsfbox{slika.8}
\end{center}
\caption{An $\ell$-face, $\ell\ge 5$, with three consecutive vertices with degrees $5,4,5$ and its reduction.}
\label{CONF_545_1}
\end{figure}

\begin{proof}
The proof follows the lines of the proof of Lemma~\ref{lm-face45x}.
We assume that $G$ has a face $f$ with vertices $v_1,\ldots,v_{\ell}$ such
that $v_1$ and $v_3$ are $5$-vertices and $v_2$ is a $4$-vertex and
we let $v'$ to be the common neighbor of $v_1$ and $v_2$ (see Figure~\ref{CONF_545_1}).
Remove the vertices $v_1$, $v_2$ and $v_3$, add the edge $v_4v_{\ell}$ and
triangulate the obtained new face with edges incident with $v'$.
The obtained graph $G'$,
which is loopless, is cyclically $(D+1)$-colorable by the minimality of $G$.

The cyclic $(D+1)$-coloring of $G'$ can now be extended to $G$.
If the color $a$ of $v'$ is not assigned to any of the vertices $v_4,\ldots,v_{\ell}$,
we color $v_3$ with $a$.
Otherwise, we color $v_3$ with any available color (as the cyclic degree of $v_3$ is $D+2$ and
two of its cyclic neighbors are uncolored, there is an available color).
We next color the vertex $v_1$ (its cyclic degree is $D+2$, it has an uncolored
cyclic neighbor and has two cyclic neighbors colored with $a$) and the vertex $v_2$ (its
cyclic degree is $D+1$ and has two cyclic neighbors colored with $a$).
\end{proof}

\section{Discharging phase}

In this section, we present the second part of the proof of our result.
At the beginning, every $d$-vertex of a $D$-minimal graph is assigned
charge of $d-6$ units and every $\ell$-face is assigned
charge of $2\ell-6$ units. The Euler formula implies that the total
amount of charge assigned to all the vertices and faces of the graph
is equal to $-12$. The initial charge is then redistributed based
on the following two rules:

\begin{description}
\item[Rule 1] 
Every $\ell$-face, $\ell \geq 4$, sends $2$ units of charge to each incident $4$-vertex. 
\item[Rule 2]
Every $\ell$-face, $\ell \geq 4$, sends $1$ unit of charge to each incident $5$-vertex.
\end{description}

First, we show that the final charge of every vertex is non-negative.

\begin{lemma}
\label{lm-vertex}
The final amount of charge of
every vertex $v$ of a $D$-minimal graph $G$, $D\ge 5$,
is non-negative.
\end{lemma}

\begin{proof}
Let $d$ be the degree of $v$. By Lemma~\ref{lm-mindeg}, the degree $d$ is at least four.
If $d\ge 6$, then the initial amount of charge of $v$ is non-negative and its final amount of charge
is also non-negative as $v$ neither receives nor sends out any charge. If $d\in\{4,5\}$,
then $v$ is incident with an $\ell$-face, $\ell\ge 4$, by Lemma~\ref{lm-deg45}.
If the degree of $v$ is four, then $v$ receive two units of charge by Rule 1, and
if its degree is five, then it receive one unit of charge by Rule 2. In either of the two
cases, the final amount of charge of $v$ is zero.
\end{proof}

We next show that the final charge of every face is non-negative.

\begin{lemma}
\label{lm-face}
The final amount of charge of
every face $f$ of a $D$-minimal graph $G$, $D\ge 5$,
is non-negative.
\end{lemma}

\begin{proof}
If $f$ is $3$-face, its final amount of charge is zero.
If $f$ is a $4$-face, then it is incident with no $4$-vertex by Lemma~\ref{lm-deg4} and
with at most two $5$-vertices by Lemma~\ref{lm-4face-555}. Hence, $f$ sends out
at most two units of charge (twice one unit by Rule 2) and the final charge of $f$
is non-negative.

Assume that $f$ is a $5$-face. By Lemma~\ref{lm-5face}, $f$ is is incident
with at most two vertices of degree four or five. Consequently, Rules 1 and 2
apply at most twice and $f$ sends out at most four units of charge to incident
vertices. Since the amount of initial charge of $f$ is equal to four units,
the final charge of $f$ is non-negative.

The remaining case is that $f$ is an $\ell$-face, $\ell\ge 6$.
Let $v_1,\ldots,v_{\ell}$ be vertices on the boundary of $f$.
If all vertices of $f$ received charge from $f$,
then they all would be $5$-vertices by Lemmas~\ref{lm-face44} and \ref{lm-face45x}.
Hence, $f$ would send out $\ell$ units of charge and its final charge
would be non-negative in this case.

In what follows, we assume that $f$ does not send out charge to all incident
vertices.
Let $A_1,\ldots,A_k$ be maximal consecutive intervals of vertices that receive
charge from $v$.
Observe that
\begin{equation}
|A_1|+\cdots+|A_k|+k\le\ell\;\mbox{.}\label{eq}
\end{equation}
We claim that the total amount of charge sent by $f$ to the vertices $v_i\in A_j$
is at most $|A_j|+1$ units for every $j=1,\ldots,k$.

If $|A_j|=1$, then $f$ sends at most two units of charge to the only vertex of $A_j$ and
the claim holds. If $|A_j|=2$, then $f$ sends at most three units of charge to the two
vertices of $A_j$ as they both cannot be $4$-vertices by Lemma~\ref{lm-face44}.
The claim also holds in this case. If $|A_j|\ge 3$, then none of the vertices
of $A_j$ is a $4$-vertex by Lemmas~\ref{lm-face44}, \ref{lm-face45x} and \ref{lm-face545}.
We conclude that $f$ sends to the vertices of $A_j$ exactly $|A_j|$ units of charge
as they all are $5$-vertices.

Since the vertices of $A_j$ receive at most $|A_j|+1$ units of charge from $f$,
the total amount of charge sent out by $f$ is at most $|A_1|+\cdots+|A_k|+k$
which is at most $\ell$. Since the initial amount of charge of $f$ is $2\ell-6$ and
$\ell\ge 6$, the final amount of charge of $f$ is non-negative.
\end{proof}

Lemmas~\ref{lm-vertex} and \ref{lm-face} yield our main result:

\begin{theorem}
\label{thm-our}
There is no $D$-minimal graph, $D\ge 5$.
Every plane graph with maximum face size $\Delta^*\ge 5$ whose all faces
of size four or more are vertex-disjoint has a cyclic coloring
with at most $\Delta^*+1$ colors.
\end{theorem}

\section*{Acknowledgement}

The authors would like to thank Riste {\v S}krekovski for discussions on cyclic colorings
of plane graphs.


\begin{thebibliography}{99}
\bibitem{bib-harman}
N. Harman:
Graphs with four independent crossings are five colorable,
manuscript.
\end{thebibliography}
\end{document}